\documentstyle[11pt]{article}
\begin{document}
\begin{center}
{\Large \bf  A bound on the number of  curves of  a given degree
through a general point of a projective variety}

\bigskip
 {\large \bf Jun-Muk Hwang}
\footnote{Supported  by the Korea Research Foundation Grant
(KRF-2002-070-C00003).


}

 \end{center}

\begin{abstract}
Let $X$ be an irreducible projective variety of dimension $n$ in a
projective space and let $x$ be a point of $X$. Denote by ${\rm
Curves}_d(X,x)$ the space of curves of degree $d$ lying on $X$ and
passing through $x$. We will show that  the number of components
of ${\rm Curves}_d(X,x)$ for any smooth point $x$ outside a
subvariety of codimension $\geq 2$ is bounded by a number
depending only on $n$ and $d$. An effective bound is given. A key
ingredient of the proof is an argument from
Ein-K\"uchle-Lazarsfeld's work on Seshadri numbers.
\end{abstract}

\bigskip

\section{Introduction}

\medskip
 This work was motivated by the following result of J.M. Landsberg's.

\medskip
{\bf Theorem 1 ([L2, Theorem 1])} {\it Let $X$ be an irreducible
projective variety of dimension $n$ in a projective space and let
$x \in X$ be a general point. Then the number of lines lying on
$X$ and passing through $x$ is either infinite or bounded by
$n!$.}

\medskip
It is remarkable that the bound $n!$ is optimal: it is  achieved
when $X$ is a smooth hypersurface of degree $n$ in ${\bf
P}_{n+1}$. However, even if we disregard the optimality of the
bound, the  uniformity of the bound is already quite remarkable.
Namely, the fact that the bound depends only on the dimension $n$
of $X$ is worth noticing. When interpreted as such a  uniform
boundedness result, Theorem 1 naturally leads to the following
questions.

\medskip
{\bf Question 1} What about curves of higher degree? Is the number
of curves of degree $d>0$ lying on $X$ and passing  through a
general point $x \in X$ either infinite or bounded by a number
depending only on $d$ and $n$?

\smallskip
{\bf Question 2} What about the case when there are infinitely
many lines through a general point $x$? Is  the number of
components of the space of lines lying on $X$ and passing through
a general point $x \in X$ bounded  by a number depending only on
$n$?

\smallskip
{\bf Question 3} What about non-general points? Is  the number of
lines lying on $X$ through any given point of $X$ either infinite
or bounded by a number depending only on $n$?

\medskip
In Landsberg's proof, the uniformity comes from his earlier result
[L1] that
  a line osculating to order $n+1$ at a general point of $X$
  must be contained in $X$.
The differential-geometric argument of [L1] using the moving frame
method seems difficult to be generalized to handle above
questions.

\medskip
 In this paper, we introduce an approach to these questions, using
tools from the study of uniform lower bounds for the Seshadri
numbers of an ample line bundle at  general points of a variety
([EKL],[HK]). Apparently, there is no direct connection between
the questions raised above and Seshadri numbers. Nonetheless  the
arguments used in [EKL] and [HK] can be translated to apply here.
Using this, we will get effective, albeit non-optimal, answers to
Question 1 and Question 2. Our result can be stated as follows.
Let us denote by ${\rm Curves}_d(X,x)$ the space of curves of
degree $d$ lying on a projective variety $X$ and passing through a
point $x \in X$ (See Notation and Conventions below for a precise
definition.)

\medskip
 {\bf Theorem 2} {\it Let $n$ and $d$ be two positive
integers. Let $X$ be an irreducible  projective variety of
dimension $n$ in a projective space and $x \in X$ be a general
point. Then the number of components of $\mbox{\rm Curves}_d(X,x)$
is bounded by
$$ \left( \begin{array}{c} (2n+2)((2nd)^n -1) \\
2n+1\end{array} \right)^{ (2n+2)(4d^2-4d+ 2)}.$$}

\medskip
Regarding Question 3, the answer is plainly no. Let us look at two
examples.

\medskip
{\bf Example 1} Let ${\bf P}_1 \times {\bf P}_2 \subset {\bf P}_5$
be the Segre embedding. Choose a curve $C$ of degree $\ell >0$ in
${\bf P}_2$. Let $X'$ be the surface ${\bf P}_1 \times C \subset
{\bf P}_1 \times {\bf P}_2$. Pick a point $o \in {\bf P}_1$ and a
point  $P$ in $\{ o\} \times {\bf P}_2 \subset {\bf P}_1 \times
{\bf P}_2$ outside $X'$. Let $X \subset {\bf P}_4$ be the
projection of $X'$ from $P$. Then the image of $\{ o\} \times C$
is a line in $X$ at a general point of which there are at least
$\ell +1$ lines. Given any dimension $n \geq 2$ and an integer
$\ell>0$, by taking the Segre product of $X$ with an irreducible
variety of dimension $n-2$ containing no lines, we get an example
of an irreducible variety of dimension $n$ where the number of
lines through a general point in a hypersurface is finite, but
larger than $\ell$.

\medskip
 {\bf Example 2} Let $k$ be an odd integer and consider the
Fermat surface $X_0^k + X_1^k + X_2^k + X_3^k =0$ in ${\bf P}_3$.
Then through the point $(1, -1, 0, 0)$ there are at least $k$
distinct lines defined by $X_0+ X_1=0$ and $X_2 + e^{\frac{2\pi
j\sqrt{-1}}{k}} X_3 =0, 1 \leq j \leq k$. Given any dimension $n
\geq 2$ and an integer $M>0$, by taking the Segre product of the
Fermat surface with a smooth variety of dimension $n-2$ containing
no lines, we get an example of a smooth variety of dimension $n$
where the number of lines through any point in a codimension 2
subset is finite, but larger than $M$.

\medskip
These examples suggest that the following result of ours gives a
more or less optimal answer to Question 3.

\medskip
{\bf Theorem 3} {\it Let $X$ ba an irreducible projective variety
of dimension $n$ in a projective space. Then there exists a
subvariety $R$ of codimension $\geq 2$ in the smooth locus of $X$
such that for any smooth point of $X$ off $R$, the number of the
components of $\mbox{\rm Curves}_d(X,x)$ is bounded by a number
$\nu_{n,d}$ depending only on $n$ and $d$. }

\medskip For an explicit value
of the bound $\nu_{n,d}$, see Definition 4.

\medskip
The rough idea of the proofs of Theorem 2 and Theorem 3 is the
following. First we will  explain in Section 2 how to obtain a
bound depending on the degree of $X$. Here the main ingredient, in
addition to some elementary projective geometry, is the effective
bound on the number of components of Chow varieties obtained in
recent works on effective bounds on the number of maps dominating
varieties of general type, e.g. [Gu] and [Ts]. Now to prove
Theorem 2 and Theorem 3,  the strategy is to study the foliation
on $X$ generated by curves of degree $d$, to be constructed in
Section 3. This foliation has the property that its general leaf
contains all curves of degree $d$ lying on $X$ passing through a
general point of the leaf, and it is the foliation of minimal rank
with this property. Thus to prove Theorem 2, we may replace $X$ by
a leaf of the foliation. The heart of the proof of Theorem 2 is to
show that the degree of the leaf can be bounded in terms of $n$
and $d$. This is achieved in Section 4 by using an argument from
Ein-K\"chle-Lazarsfeld's work on Seshadri numbers ([EKL]). The
proof of Theorem 3, presented in Section 5,  is  by an induction
argument using Theorem 2 and by a study of the foliation in
codimension 1.

\bigskip
{\bf Notation and Conventions}

\medskip
1.  Throughout this paper, we will work over the complex numbers.

2. A variety need not be irreducible,  but has finitely many
components.  For a variety $Y$, $\#Y$ denotes the number of
irreducible components of $Y$.

3. By a curve of degree $d$ on an irreducible projective variety
$Y$ in a projective space, we mean an irreducible reduced
subvariety of dimension 1 lying on $Y$ which has degree $d$ with
respect to the hyperplane line bundle of the projective space.

4. For an irreducible subvariety $Y$ in a projective space, denote
by ${\rm Chow}_{1,d}(Y)$ the Chow variety  of effective 1-cycles
of degree $d$ on $Y$ (see [Ko, I.3] for the definition). Denote by
${\rm Curves}_d(Y)$ the quasi-projective subvariety   in ${\rm
Chow}_{1,d}(Y)$ parametrizing curves of degree $d$ on $Y$. For a
point $y \in Y$, denote by ${\rm Curves}_d(Y,y)$ the subvariety of
${\rm Curves}_d(Y)$ parametrizing members passing through $y$.

5. We will say that a property holds at a general point of an
irreducible quasi-projective variety $Y$ if it holds for a
non-empty Zariski-open subset of $Y$. It holds at a very general
point  if it is satisfied off the union of countably many proper
closed subvarieties of $Y$.

\section{A bound depending on the degree of the variety }

Let us start by recalling the following elementary fact.

\medskip
{\bf Proposition 1} {\it Let $X \subset {\bf P}_N$ be an
irreducible projective variety and let $x \in X$ be a smooth
point. Let $\pi: \tilde{X} \rightarrow X$ be the blow-up of $X$ at
$x$ and let $E$ be the exceptional divisor. Denote by $H$  the
hyperplane  divisor on $ {\bf P}_N$. Then $2 \pi^* H -E$ is a very
ample divisor on $\tilde{X}$. If the degree of $X$ is $a$, the
degree of $\tilde{X}$ with respect to $2 \pi^*H-E$ is $2^n a -1$.}

\medskip
{\it Proof}. Let $\pi': \tilde{\bf P}_N \rightarrow {\bf P}_N$ be
the blow-up of ${\bf P}_N$ at $x$ and let $E' \subset \tilde{\bf
P}_N$ be the exceptional divisor. More precisely, $\tilde{\bf
P}_N$ is the subvariety of ${\bf P}_N \times {\bf P}_{N-1}$ which
is the closure of the graph of the projection of ${\bf P}_N$ to a
hyperplane ${\bf P}_{N-1}$ with the vertex $x$ and $E'$ is the
inverse image of $x$ in $\tilde{\bf P}_N$ under the projection
$p_1: \tilde{\bf P}_N \rightarrow {\bf P}_N$. Then $E'$ is
biregular to ${\bf P}_{N-1}$ by the projection $p_2: \tilde{\bf
P}_N \rightarrow {\bf P}_{N-1}$. Let $H$ be the hyperplane divisor
on ${\bf P}_N$ and let $H'$ be the hyperplane divisor on ${\bf
P}_{N-1}$. Then the divisor
$$(p_1^*H + p_2^*H')|_{\tilde{\bf P}_N} = 2 \pi'^* H -E'$$ is
very ample on $\tilde{\bf P}_N$. The proper image of $X$ in
$\tilde{\bf P}_N$ is $\tilde{X}$ and the divisor $E \subset
\tilde{X}$ is just the restriction of $E'$ to $\tilde{X}$.
 Thus
$2 \pi^*H-E =(2 \pi'^*H-E')|_{\tilde{X}}$ is very ample on
$\tilde{X}$. The degree is
$$ (2 \pi^*H -E)^n = 2^n a -1$$ from $a =
H^n$ and $E^n = (-1)^{n-1}$. $\Box$

\medskip
{\bf Proposition 2} {\it Let $x$ be a smooth point of an
irreducible projective variety $X \subset {\bf P}_N$. Assume that
a general member of each component of $\mbox{\rm Curves}_d(X,x)$
has multiplicity $\leq m$ at $x$ for some positive integer $m$.
Let $\pi: \tilde{X} \rightarrow X$ be the blow-up of $X$ at $x$.
Identify $\tilde{X}$ as a subvariety of a projective space via the
very ample divisor $2 \pi^*H-E$ of Proposition 1. Then
$$\#\mbox{\rm Curves}_d(X,x) \leq \sum_{i=1}^{m} \#\mbox{\rm Curves}_{2d-i}(\tilde{X}).$$}

\medskip
{\it Proof}. Fix an integer $i, 1 \leq i \leq m$. Let ${\rm
Curves}'_{2d-i}(\tilde{X})$ be the union of the components of
${\rm Curves}_{2d-i}(\tilde{X})$ whose members have intersection
number $i$ with $E$. Any member of ${\rm
Curves}'_{2d-i}(\tilde{X})$ is sent by $\pi$ to a curve of degree
$d$ passing through $x$ which has multiplicity $i$ at $x$. This
induces a morphism $$\pi_*: \bigcup_{i=1}^m {\rm
Curves}'_{2d-i}(\tilde{X}) \rightarrow {\rm Curves}_d(X,x).$$
Since a general member of each irreducible component of ${\rm
Curves}_d(X,x)$ has multiplicity $\leq m$, the morphism $\pi_*$ is
dominant on each irreducible component of ${\rm Curves}_d(X,x)$.
Thus $$\#{\rm Curves}_d(X,x) \leq \sum_{i=1}^m \#{\rm
Curves}'_{2d-i}(\tilde{X}) \leq \sum_{i=1}^m \#{\rm
Curves}_{2d-i}(\tilde{X}). \;\;\; \Box$$

\medskip
{\bf Proposition 3} {\it Let $X \subset {\bf P}_N$ be an
irreducible projective variety and let $x$ be a general point of
$X$. Then a general member of any component of ${\rm
Curves}_d(X,x)$  is smooth at $x$.}

\medskip
{\it Proof}. Let $B$ be a component of ${\rm Curves}_d(X)$ whose
members sweep out an open subset of $X$ and let $G \subset B
\times X$ be the subvariety defined by the incidence relation
$$G = \{ (b, x) \in B \times X, x \in C_b \}$$
where $C_b$ denotes the curve in $X$ corresponding to $b \in B$.
It is clear that $G$ contains a unique irreducible component $G^o$
of dimension $\dim B +1$.  There exists an open subset $X' \subset
X$ such that for each $x \in X'$, the intersection $pr_X^{-1}(x)
\cap G$ is contained in $G^o$. By shrinking $X'$, we may assume
that each component of $pr_X^{-1}(x) \cap G $ has dimension $\dim
B + 1 - \dim X$ for each $x \in X'$. Let $S_x \subset G \cap
pr_X^{-1}(x)$ be the subvariety defined by
$$S_x := \{(b, x), C_b \mbox{ is singular at } x \}.$$ Suppose
$S_x$ contains some components of $G \cap pr_X^{-1}(x)$ for
general $x \in X'$. Then the union of $S_x$'s as $x$ varies define
a subset of $G$ of dimension $\dim B +1$. This implies that for a
general $x \in X'$, we have the equality $S_x = G \cap
pr_X^{-1}(x)$. This is absurd because some member of $B$ must be
smooth at $x$. Hence $S_x$ contains no component of $G \cap
pr_X^{-1}(x)$ for a general $x \in X$.  Since this is true for any
finitely many possible choices of the component $B$, Proposition 3
is proved. $\Box$

\medskip
We will use an effective bound on the number of components of the
Chow variety. There are a number of results  obtained in this
direction, e.g. [Gu] and [Ts]. For example, the bound given in
[Gu] implies  the following.

\medskip
{\bf Proposition 4 ([He, Proposition 3.6])} {\it The number of
components of the Chow variety ${\rm Chow}_{k, \delta}(Y)$
parametrizing subvarieties of dimension $k$ and degree $\delta$ in
a projective variety $Y$ of degree $\delta'$ in ${\bf P}_{N'}$ is
bounded by
$$\left( \begin{array}{c} (N'+1) \max \{\delta, \delta'\} \\
N'\end{array} \right)^{ (N'+1) \left[\delta \left(
\begin{array}{c} \delta+ k-1 \\k \end{array} \right) + \left(
\begin{array}{c} \delta+ k -1 \\k-1 \end{array} \right)
\right]}.$$}

\medskip
{\bf Definition 1}  Given positive integers $n, d$ and $a$,
 define \begin{eqnarray*}
\mu_{n,d,a} &:=& \sum_{i=1}^d \left(
\begin{array}{c} (2n+2)
 \max \{ 2d-i, 2^n a -1 \} \\
2n+1\end{array} \right)^{ (2n+2)(4d^2-4di+ i^2 + 1)} \\
\lambda_{n,d,a} &:=& \left(
\begin{array}{c} (2n+2)
 \max \{ 2d-1, 2^n a -1 \} \\
2n+1\end{array} \right)^{ (2n+2)(4d^2-4d+ 2)}. \end{eqnarray*}

\medskip
 {\bf Proposition 5} {\it
 Let $X \subset {\bf P}_N$ be an
irreducible projective variety of dimension $n$ and degree $a$.
Let $x$ be a smooth point of $X$. Then $$\# \mbox{\rm
Curves}_d(X,x) \leq \mu_{n,d,a}.$$ If furthermore $x$ is a general
point of $X$, then
$$ \#\mbox{\rm Curves}_d(X,x) \leq \lambda_{n,d,a}.$$}

\medskip
{\it Proof}. Using Proposition 2, we will bound $\# {\rm
Curves}_{2d-i}(\tilde{X})$ for $1\leq i \leq m$. We may project
$\tilde{X}$ to ${\bf P}_{2n+1}$ to count $\#{\rm
Curves}_{2d-i}(\tilde{X})$. Thus we can use Proposition 4 with
$N'= 2n+1, \delta=2d-i, \delta'= 2^n a -1$ and $k=1$.
 Noting that the multiplicity at $x$ of
any member of ${\rm Curves}_d(X,x)$ is bounded by $m=d$, we get
the first inequality.  For a general point $x \in X$, we can set
$m=1$ in Proposition 2 by Proposition 3, which gives the second
inequality. $\Box$

\section{Foliation generated by curves of degree $d$}

\medskip
Let $X$ be an irreducible projective variety of dimension $n$ in a
projective space. Fix a positive integer $d$. The goal of this
section is to construct a foliation of minimal rank on an open
subset of $X$ such that members of ${\rm Curves}_d(X,x)$ lies in
the closure of the leaf through $x$ for a general $x \in X$. The
construction is similar to the construction of
Seshadri-exceptional foliation in [HK]. Here we will do it with
some more care because we have to study the foliation in
codimension 1 for Theorem 3. Let us start with the definition of a
foliation. This may not be the standard definition, but it will be
convenient for us.

\medskip
 {\bf Definition 2} Let $X$ be an irreducible projective
variety. Denote by $Sm(X)$ the smooth locus of $X$. A subsheaf
${\cal F}$ of the tangent sheaf $T_{Sm(X)}$ of $Sm(X)$ is called a
foliation on $X$ if it satisfies the following two conditions.

 (1) The quotient
 $T_{Sm(X)}/{\cal F}$  is torsion-free on $Sm(X)$. This implies that the open
 subset $${\rm Dom}({\cal F}) := \{ x \in Sm(X), T_{Sm(X)}/{\cal F} \mbox{ is locally free
 at } x  \}$$   is the complement of
 a subvariety of codimension $\geq 2 $ in $Sm(X)$.

 (2)  For each $x \in {\rm
 Dom}({\cal F})$,
there exists a complex analytic submanifold ${\cal F}_x \subset
{\rm Dom}({\cal F})$ passing through $x$, called the $x$-leaf of
${\cal F}$, such that the fiber of ${\cal F}$ at each point $y$ of
${\cal F}_x$ is the tangent space of ${\cal F}_x$ at $y$.

\medskip
{\bf Proposition 6} {\it Let $X$ be an irreducible projective
variety. Suppose for each very general point $x \in X$, we can
assign an irreducible projective subvariety $Z_x$ such that for a
very general point $x \in X$ and  a very general point $y \in
Z_x$, the two subvarieties $Z_x$ and $Z_y$ coincide. Then there
exists a unique foliation ${\cal F}$ on $X$ such that the $x$-leaf
${\cal F}_x$ is an open subset of $Z_x$ for a very general $x$.}

\medskip
{\it Proof}. By the countability of the components of the Hilbert
scheme of $X$, there exists an irreducible family of subvarieties
$\rho: {\cal U} \rightarrow {\cal D}$ for some sub-scheme ${\cal
D} \subset {\rm Hilb}(X)$ such that for a very general point $x
\in X$, the subvariety $Z_x$ is the image of a fiber of $\rho$ by
the evaluation morphism $\eta: {\cal U} \rightarrow X$. The
assumption that $Z_x =Z_y$ for a very general $x \in X$ and a very
general $y \in Z_x$ implies that $\eta$ is birational. Consider
the subsheaf ${\cal F}'$ of the tangent sheaf of $Sm(X)$ defined
by the push-forward of the relative tangent sheaf of $\rho$
restricted to $\eta^{-1}(Sm(X))$. Define ${\cal F}$ to be the
double-dual of ${\cal F}'$. It is obvious that ${\cal F}$ is the
unique foliation satisfying the desired property. $\Box$

\medskip
 {\bf Proposition 7} {\it Let
$X$ be an irreducible projective variety of dimension $n$ in a
projective space. Fix a positive integer $d$. Then there exists a
unique foliation ${\cal F}$ on $X$, which we will call the
foliation generated by curves of degree $d$, with the following
properties.

 (i) Each leaf of the foliation is a smooth quasi-projective
subvariety in ${\rm Dom}({\cal F})$.

(ii) For a general point $x \in X$, all members of $\mbox{\rm
Curves}_d(X,x)$ are contained in the closure of the leaf of the
foliation passing through $x$.

(iii) The foliation ${\cal F}$ is minimal with respect to (ii). In
other words, if there exists a  foliation ${\cal G}$ satisfying
(ii) then  ${\cal F}_x \subset {\cal G}_x$ for a general point $x
\in X$.}

\medskip
{\it Proof}.   For each $x \in X$, define
$$C_x := \mbox{closure of } \bigcup_{[C] \in {\rm Curves}_d(X,x)} C.$$
Then $C_x$ is a projective subvariety in $X$, not necessarily of
pure dimension. Choose a smooth affine open subset $T \subset X$
and consider the incidence relation $ {\cal C} \subset T \times X$
defined by
$${\cal C} = \{(t, x), x \in C_t \}.$$ We may assume that the projection
$pr_T: {\cal C} \rightarrow T$ is flat by replacing $T$ by an open
subset of $T$. Let us make this assumption. For an irreducible
subvariety $W \subset X$ intersecting $T$, let
$$C_W :=\mbox{ closure of } \bigcup_{s \in W \cap T} C_s.$$  This
subvariety $C_W$ of $X$ is not necessarily irreducible. But each
component of $C_W$ contains $W$ by the flatness of $pr_T: {\cal C}
\rightarrow T $. This implies that either every component of $C_W$
has dimension strictly larger than $W$ or $C_W=W$. Note that when
$W$ is one point $x \in T$, $C_W = C_x$.

 An irreducible subvariety $W$ is said to be
saturated if $ C_W =W$.
  For each $x \in T$, there exists a unique minimal saturated
subvariety $Z_x$ containing $x$ constructed as follows. Let
$Z^1_x$ be a component of $C_x $ and inductively define
$Z^{i+1}_x$ to be a component of $C_{Z^i_x}$. Then
$\dim(Z^{i+1}_x) > \dim (Z^i_x)$ or $Z^i_x = Z^{i+1}_x=Z^{i+2}_x =
\cdots$. Thus $Z^n_x = Z^{n+1}_x =\cdots$. Define $Z_x = Z^n_x$.
Then $Z_x$ is saturated. We claim that any saturated subvariety
containing $x$ contains $Z_x$. In fact, if $W$ is a saturated
subvariety and $W' \subset W$ is any irreducible subvariety of $W$
intersecting $T$, $C_{W'} \subset W$. Thus if $x \in W$, $Z^i_x
\subset W$ inductively for all $i$.

We claim that $Z_x =Z_y$ for a very general $x \in X$ and a
 very general $y \in Z_x$. Assuming the claim, let us finish the
 proof of Proposition 7. By Proposition 6, we get a foliation
 ${\cal F}$ on $X$ such that the $x$-leaf  ${\cal F}_x$ is an open
 subset of $Z_x$ for a very general $x$. This implies that for a very general $x$,
  the $x$-leaf  ${\cal F}_x$ is
quasi-projective and
 contains all members of ${\rm Curves}_d(X,x)$. But then
  the same holds for a general $x$.  This shows (i) and  (ii).
 The condition (iii) follows from
 the fact that $Z_x$ is the minimal saturated subvariety containing
 $x$ for a general $x$.

Now to prove the claim, notice that $\dim Z_x = \dim Z_y$ for very
general $x$ and $y$. For a very general $x \in X$ and a very
general $y \in Z_x$,  we have $Z_y \subset Z_x$ by the
saturatedness  of $Z_x$. Thus $Z_x=Z_y$ by $\dim Z_x = \dim Z_y$.
$\Box$

\section{Proof of Theorem 2}

\smallskip
{\bf Definition 3} Given positive integers $n$ and $d$, define
$$\lambda_{n,d} := \lambda_{n,d, (nd)^n} = \left( \begin{array}{c} (2n+2)((2nd)^n -1) \\
2n+1\end{array} \right)^{ (2n+2)(4d^2-4d+ 2)}.$$ For a fixed $d$,
$\lambda_{n,d}$ is an increasing function of $n$.

\medskip
We can restate  Theorem 2 as follows.

\medskip
 {\bf Theorem 2} {\it Let $n$ and $d$ be two positive
integers. Let $X$ be an irreducible projective variety of
dimension $n$ in a projective space and $x \in X$ be a general
point. Then
$$\# \mbox{\rm Curves}_d(X,x) \leq \lambda_{n,d}.$$ }

\medskip
The heart of the proof of Theorem 2 is the following.

\medskip
{\bf Proposition 8} {\it Let $X$ be a projective variety of
dimension $n$ in ${\bf P}_N$. Fix a positive integer $d$. Assume
that the foliation ${\cal F}$ generated by curves of degree $d$
has rank $n$.  Then the degree of $X$ is bounded by $(nd)^n$.}

\medskip
{\it Proof}. The proof is a translation of [HK, Proof of Theorem
1].  We start by recalling the definition of the multiplicity of
an effective divisor along an irreducible subvariety. For an
effective divisor $D$ on a variety $Y$ and a non-singular point $y
\in Y$, let $m_y(D)$ denote the multiplicity of $D$ at $y$. For an
irreducible subvariety $Z \subset Y$ intersecting the smooth locus
of $Y$, let $m_Z(D)$ denote $m_y(D)$ at a general point $y \in Z$.

Now assume that the degree of $X$ is strictly bigger than
$(nd)^n$. Using the notation in the proof of Proposition 7,  let
$\Gamma \subset X \times T$ be the closure of the graph of the
inclusion $T \subset X$. {}From the relation between the degree
and the Hilbert polynomial of $X$,
$$\dim H^0(X, {\cal O}(kH)) \geq \frac{(nd)^n}{n!} k^n +
O(k^{n-1}).$$ As in [EKL, (3.8)], this implies that
 if $k \gg 0$, there exists a
divisor $D \in |{\cal O}_{X \times T}(pr_1^*(kH))|$ with
$m_{\Gamma}(E) > kdn$.

The following lemma is  essentially equal to [EKL, Lemma 3.5.1] or
[HK, Lemma 1]. Its proof will be omitted.

\medskip
{\bf Lemma 1} {\it Let $Z \subset X \times T$ be an irreducible
closed subvariety dominating both $X$ and $T$. Then there exists
an irreducible closed subvariety $CZ \subset X \times T$ having
the following properties:

(i) $Z \subset CZ$.

(ii) For general $t \in T$ with the fiber $Z_t$ intersecting $T$,
the fiber $(CZ)_t \subset X$ consists of some components of $
C_{Z_t}$. Here $Z_t$ is regarded as a subvariety of $X$ by the
projection $X \times T \rightarrow X$ and $C_{Z_t}$ denotes the
subvariety constructed in the proof of Proposition 7.}

\medskip
Using Lemma 1 as in [EKL, (3.7)], we construct a nested sequence
of
 irreducible  subvarieties $$Z_0 \subset Z_1 \subset \cdots \subset Z_i
 \subset \cdots \subset Z_n$$ in $  X\times T$ as follows.
 Set $Z_0= \Gamma$. Then inductively
define $Z_{i+1} = CZ_{i}$. {}From the construction, $(Z_n)_t =
Z_{t}$ for a general $t \in T$ and $Z_{n+1}= Z_{n}$. By the
assumption that the rank of ${\cal F}$ is $n$, we have $X=Z_x$ for
a very general $x$. This implies that $Z_n= T \times X$.

Consider the multiplicities $m_{Z_i}(D)$. We have
$$m_{Z_0}(D) = m_{\Gamma}(D) >
 kdn,\; \;  m_{Z_n}(D) = m_{X\times T}(D)=0.$$ It follows that there is
 at least one index $i, 0 \leq i \leq n-1$ such that $$m_{Z_i}(D)
 - m_{Z_{i+1}}(D) > kd.$$

Now we  use the following result from [EKL].

\medskip
{\bf Lemma 2 ([EKL, Proposition 2.3.])} {\it Let $Y$ and $T$ be
smooth irreducible  varieties, with $T$ affine, and suppose that
$Z \subset V \subset Y \times T$ are irreducible subvarieties such
that $V$ dominates $Y$. Let $L$ be a line bundle on $Y$, and
suppose given on $Y \times T$ a divisor $D \in |pr_Y^*L|$. Then
there exists a divisor $D' \in |pr_Y^*L|$ on $Y \times T$ whose
support does not contain $V$ such that $$ m_Z(D') \geq m_Z(D) -
m_V(D).$$}

\medskip
First, let us assume that $X$ is smooth. By Lemma 2, applied to $L
=H$ and $Y=X$, there exists a divisor $D' \in |{\cal O}_{X \times
T}(pr_X^*(kH))|$ such that $ m_{Z_i}(D') > kd$ and $ Z_{i+1}$ is
not contained in the support of $D'$. Then for a general $t \in T$
and $x \in (Z_i)_t$, $m_x(D'_t)
> kd$ as in [EKL,(3.9)]. But there exists a  curve $C'$ of degree $d$ on $X$
 passing through
 $x$ which is contained in $(Z_{i+1})_t$, but not contained in
the support of $D'_t$. Thus we get the contradiction
$$ kd = k (H \cdot C') = D'_t \cdot C' \geq m_x(D'_t)  >
kd. $$ Now when $X$ has singularity, apply the same argument to
$Y= \hat{X}$, a desingularization of $X$. We have
$$m_{\hat{Z}_i}(\hat{D}) - m_{\hat{Z}_{i+1}}(\hat{D}) >kd$$
for the proper images $\hat{Z}_i$ and $\hat{Z}_{i+1}$ in $\hat{X}$
of $Z_i$ and $Z_{i+1}$ and the pull-back divisor $\hat{D}$. Thus
the same contradiction occurs. $\Box$

\medskip
{\it Proof of Theorem 2}.  We will use induction on the dimension
$n$. It is obvious for $n = 1$. Suppose the rank of ${\cal F}$ is
$n$.  Then the bound follows from  Proposition 5 with $a = (nd)^n$
by Proposition 8. Suppose the rank of ${\cal F} $ is $r <n$. If
$r=0$, or equivalently, ${\rm Curves}_d(X,x)$ is empty for a
general $x$, there is nothing to prove. So let us assume that $r
\geq 1$. For a general $x \in X$, we may assume that $x$ is a
general point of ${\cal F}_y$ for some $y \in X$. By Proposition 7
(ii)
$$\#{\rm Curves}_d(X,x) = \# {\rm Curves}_d({\cal F}_y,x).$$ But
 $$\# {\rm Curves}_d({\cal F}_y,x) \leq \lambda_{r,d} \leq \lambda_{n,d}$$ by the
 induction hypothesis.  This proves Theorem 2. $\Box$

\section{Proof of Theorem 3}

\smallskip
{\bf Definition 4} Given positive integers $n$ and $d$, define
inductively \begin{eqnarray*} \nu_{1,d} &:=& 1 \\ \nu_{n,d} &:=&
\max\{ \lambda_{n-1,d} + \nu_{n-1,d}, \;\; \mu_{n,d,(nd)^n}
\}.\end{eqnarray*} For a fixed $d$, it is an increasing function
of $n$.

\medskip
We can restate Theorem 3 as follows.

\medskip
{\bf Theorem 3} {\it Let $X$ be an irreducible projective variety
of dimension $n$ in a projective space. Let $R$ be the subvariety
defined by $$R := \{ x\in Sm(X), \# {\rm Curves}_d(X,x) >\nu_{n,d}
\}.$$ Then  the codimension of $ R \subset Sm(X)$ is $\geq 2$.}

\medskip
{\it Proof}. We will use induction on the dimension $n$. Theorem 3
is obvious for $n = 1$.

 Let ${\cal F}$ be the foliation generated by
the curves of degree $d$. If the rank of ${\cal F}$ is $n$, the
degree of $X$ is bounded by $(nd)^n$. Thus by Proposition 5, we
get
$$ \# {\rm Curves}_d(X,x) \leq \mu_{n,d, (nd)^n}$$ at any smooth point $x \in X$.
In other words, $R = \emptyset$.

 Now suppose the rank of ${\cal
F}$ is $r <m$ and  $R$ contains a component $V$ of codimension 1.
Let $v $ be a general point of $V$. Then $v \in {\rm Dom }({\cal
F})$ from Definition 2 (1). Let $F$ be the closure of the $v$-leaf
${\cal F}_v$. When $r
>0$, the set $F \cap V$ is of codimension $\leq 1$ in $F$. By
slightly moving $v$, we may assume that $v$ is a smooth point of
$F$ outside any given subset of codimension $\geq 2$ in $F$. Then
by the induction hypothesis,
$$ \#{\rm Curves}_d(F,v) \leq \nu_{r,d}.$$ The same inequality holds
for the case $r=0$ because $ \#{\rm Curves}_d(F,v) =0$ in that
case.

Let ${\rm Curves}^s_d(X)$ be the union of components of ${\rm
Curves}_d(X)$ whose members sweep out an open subset of $X$, the
superscript $s$ indicating `sweeping'. Let ${\rm
Curves}^{ns}_d(X)$ be the union of the other components of ${\rm
Curves}_d(X)$, the superscript $ns$ indicating `non-sweeping'. For
a point $x \in X$, define
$${\rm Curves}^s_d(X,x) := {\rm Curves}_d(X,x) \cap {\rm Curves}_d^s(X)$$
$${\rm Curves}^{ns}_d(X,x) := {\rm Curves}_d(X,x) \cap {\rm
Curves}_d^{ns}(X).$$ By definition, $$ {\rm Curves}_d(X,x) = {\rm
Curves}^s_d(X,x) \cup {\rm Curves}_d^{ns}(X,x)$$ which is not
necessarily a disjoint union.

\medskip
{\bf Lemma 3} {\it The subvariety $$\{ x \in X, \#{\rm
Curves}_d^{ns}(X,x)
> \lambda_{n-1, d} \}$$ has codimension $\geq 2 $ in $X$. }

\medskip
{\it Proof}. Let ${\rm Loc}^{ns}$ be the subvariety of dimension
$\leq n-1$ in $X$ which is the closure of the  union  of members
of ${\rm Curves}^{ns}_d(X)$. If $y \in X - {\rm Loc}^{ns}$, then
$${\rm Curves}_d^{ns}(X,y) = \emptyset.$$ At a general point $y$
of a component $Y$ of ${\rm Loc}^{ns}$, $$\#{\rm Curves}_d(Y, y)
\leq \lambda_{n-1,d}$$ by Theorem 2. Since there is no other
component of ${\rm Loc}^{ns}$ containing $y$ by the generality of
$y$,
$$\#{\rm Curves}_d^{ns}(X, y) = \#{\rm Curves}_d(Y,y) \leq
\lambda_{n-1,d}.$$ This proves the lemma. $\Box$

\medskip
{\bf Lemma 4} {\it Suppose $x$ is a point on $\mbox{\rm Dom}({\cal
F})$, then any member of ${\rm Curves}_d^s(X,x)$ is contained in
the closure $F$ of the $x$-leaf ${\cal F}_x$, in other words,
$${\rm Curves}_d^s(X,x) \subset {\rm Curves}_d(F,x).$$ In particular,
$\#{\rm Curves}_d^s(X,x) \leq \#{\rm Curves}_d(F,x)$ if $x \in
{\rm Dom}({\cal F})$.}

\medskip
{\it Proof}. Suppose not. Let $C$ be a member of ${\rm
Curves}_d^s(X,x)$ which is not contained in $F$. Deformations of
$C$ sweep out an open subset of $X$ by the definition of ${\rm
Curves}_d^s(X,x)$. Since $C$ is not contained in $F$, a general
deformation of $C$ will not be tangent to the foliation ${\cal
F}$. Thus at a general point $y \in X$, there is a curve of degree
$d$ through $y$ which is not contained in the closure of the leaf
of ${\cal F}$ through $y$. This is a contradiction to Proposition
7 (ii). $\Box$

\medskip
{}From Lemma 3 and Lemma 4,  if $v$ is a general point of $ V\cap
{\rm Dom}({\cal F})$, then
\begin{eqnarray*} \#{\rm Curves}_d(X,v) &\leq&
\#{\rm Curves}^{ns}_d(X,v) + \#{\rm Curves}^s_d(X,v) \\
&=& \#{\rm Curves}^{ns}_d(X,v) + \# {\rm Curves}_d(F,v)
\\ &\leq & \lambda_{n-1, d} + \nu_{r, d}\\ &\leq & \lambda_{n-1,d} + \nu_{n-1,d} \\
 &\leq& \nu_{n,d}.
\end{eqnarray*}
This contradicts the choice of $v \in R$. $\Box$

\bigskip
{\bf Acknowledgment} This work was done while the author was
visiting Mathematics Department of Harvard University for the year
2003-2004. He is grateful to the hospitality and the financial
support provided by Professor Yum-Tong Siu and Harvard University.
He would like to thank Dr. Gordon Heier for showing him the
preprint [He] and the references [Gu], [Ts].

\bigskip
{\bf References}

\medskip

[EKL] Ein, L., K\"uchle, O., Lazarsfeld, R.: Local positivity of
ample line bundles. J. Diff. Geom. {\bf 42} (1995) 193-219

[Gu] Guerra, L.: Complexity of Chow varieties and number of
morphisms on surfaces of general type. Manuscripta Math. {\bf 98}
(1999) 1-8

[He] Heier, G.: Effective finiteness theorems for maps between
canonically polarized compact complex manifolds. preprint.
alg-geom/0311086

 [HK] Hwang,
J.-M., Keum, J.: Seshadri-exceptional foliations. Math. Annalen
{\bf 325} (2003) 287-297

 [Ko] Koll\'ar, J.: {\it Rational curves on algebraic
varieties.} Ergebnisse der Mathematik und ihrer Grenzgebiete, 3
Folge, Band 32, Springer Verlag, 1996

 [L1] Landsberg, J. M.: Is a linear space contained in a
 submanifold?--On the number of derivatives needed to tell. J. reine angew.
  Math. {\bf 508} (1999) 53-60

[L2] Landsberg, J. M.: Lines on projective varieties. J. reine
angew. Math. {\bf 562} (2003) 1-3

[Ts] Tsai, I.-H.: Chow varieties and finiteness theorems for
dominant maps. J. Algebraic Geom. {\bf 7} (1998) 611-625

\bigskip
Korea Institute for Advanced Study

207-43 Cheongryangri-dong

Seoul, 130-722,  Korea

jmhwang@kias.re.kr
\end{document}